\documentclass[a4paper,11pt]{amsart}
\usepackage{fancyhdr}
\usepackage{amsmath}
\usepackage{dsfont}
\usepackage{hyperref}
\usepackage[mathscr]{eucal}
\usepackage[cp1251]{inputenc}
\usepackage[english]{babel}
\usepackage{enumerate,float,indentfirst}
\usepackage{graphicx}
\usepackage{xcolor}
\usepackage{latexsym,a4,mathrsfs,amsthm,amsmath,amssymb,url}
\usepackage{amsfonts}
\usepackage{amssymb}
\usepackage{caption}
\usepackage{blkarray}

\usepackage{booktabs}

\usepackage{stackengine}

\numberwithin{equation}{section}
\newtheorem{theorem}{Theorem}[section]

\newtheorem{proposition}{Proposition}[section]
\newtheorem{definition}{Definition}[section]
\newtheorem*{hypothesis}{Hypothesis}
\newtheorem*{problem}{Problem}
\newtheorem{remark}{Remark}[section]
\newtheorem{lemma}{Lemma}[section]

\newcommand{\R}{\mathbb R}

\newcommand{\col}{\operatorname{col}}
\newcommand{\diag}{\operatorname{diag}}

\newcommand{\sgn}{\operatorname{sgn}}
\newcommand{\eps}{\varepsilon}
\newcommand{\Energy}{\mathcal E}
\newcommand{\SpTrees}{\mathcal T}
\newcommand{\chain}{\mathcal C}

\begin{document}

\title[About subspaces the most deviating from the coordinate ones]{\bf About subspaces the most deviating \\ from the coordinate ones }

\author[Yu. Nesterenko]{Yuri Nesterenko}
\email{yuri.r.nesterenko@gmail.com}

\begin{abstract}
Using the largest principal angle as a distance between same-dimensional linear subspaces of $\mathbb{R}^n$, we construct $k$-dimensional subspaces which deviate from every coordinate $k$-subspace by at least $\arccos(1/\sqrt n)$. The construction is motivated by the hypothesis of Goreinov, Tyrtyshnikov and Zamarashkin that this value is the largest possible one for all $n > k > 0$. The subspaces are scaled star spaces of $2$-connected series--parallel graphs with $k+1$ vertices and $n$ edges, equipped with a particular positive edge weighting, while the largest principal angles take two values -- $\arccos(1 / \sqrt{n})$ and $\pi/2$, depending on whether a $k$-edge subgraph corresponding to a coordinate $k$-subspace is a spanning tree or not.

For a fixed series--parallel graph, we also prove that the constructed weighting is the unique positive one, up to scaling, for which the corresponding $k$-subspace deviates from all coordinate $k$-subspaces by at least $\arccos(1 / \sqrt{n})$.
\end{abstract}

\maketitle

\thispagestyle{empty}

\section{Introduction}
We study the real case of the problem of finding linear subspaces that deviate as much as possible from all coordinate subspaces of the same dimension.

\begin{problem}
Given $n>k>0$, what is the largest possible value of
\[
   \min_{E_I}\theta_{\max}(U,E_I),
\]
where $U$ is a $k$-dimensional subspace of $\R^n$, the minimum is taken over all coordinate $k$-subspaces $E_I \subset \R^n$, and $\theta_{\max}$ denotes the largest principal angle?
\end{problem}

The problem is closely related to a hypothesis formulated in \cite{GTZ1997} about submatrices with the best-bounded inverses.

\begin{hypothesis}
For every $n > k > 0$ and every real $n \times k$ matrix $A$ with orthonormal columns, there exists a $k \times k$ submatrix of $A$ whose inverse has spectral norm at most $\sqrt{n}$.
\end{hypothesis}

Equivalently, for every such $A$ there is a set $I\subset\{1,\ldots,n\}$ with
$|I|=k$ such that
\[
   \sigma_{\min}(A_I)\ge \frac1{\sqrt n}.
\]
If the columns of $A$ span $U$ and $E_I$ is the corresponding coordinate subspace, the singular values of $A_I$ are the cosines of the principal angles between $U$ and $E_I$. Thus the hypothesis says that the answer to the problem is at most $\arccos(1/\sqrt n)$.

The purpose of this paper is to describe a finite class of subspaces attaining this value for every $n > k > 0$, and to prove that this value cannot be increased within a certain wider class of subspaces.

The construction arose from numerical experiments. For each $1 \le k < n \le 9$ we optimized, over all $k$-dimensional subspaces of $\R^n$, the target function
\begin{equation}\label{tf}
   U \longmapsto \min_{E_I}\theta_{\max}(U,E_I).
\end{equation}
All observed maxima were equal to $\arccos(1/\sqrt n)$. Modulo coordinate permutations and sign changes, the number of observed symmetry classes was as follows.

\begin{table}[h]
\begin{tabular}{c|cccccccc}
$n\backslash k$ & $1$ & $2$ & $3$ & $4$ & $5$ & $6$ & $7$ & $8$ \\
\midrule
$2$ & $1$ \\
$3$ & $1$ & $1$ \\
$4$ & $1$ & $1$ & $1$ \\
$5$ & $1$ & $2$ & $2$ & $1$ \\
$6$ & $1$ & $3$ & $4$ & $3$ & $1$ \\
$7$ & $1$ & $4$ & $8$ & $8$ & $4$ & $1$ \\
$8$ & $1$ & $5$ & $14$ & $19$ & $14$ & $5$ & $1$ \\
$9$ & $1$ & $6$ & $23$ & $42$ & $42$ & $23$ & $6$ & $1$ \\
\end{tabular}
\caption{Numbers of symmetry classes of subspaces attaining $\arccos(1/\sqrt n)$ in numerical experiments for various $n$ and $k$.}
\label{tab}
\end{table}

Searching for this triangular array in the On-Line Encyclopedia of Integer Sequences \cite{OEIS} yielded one exact match: "Triangle read by rows: number of isomorphism classes of series--parallel matroids of rank d on n elements" \cite{A115594}, where $d = k$ in our notation.

Let us give several definitions needed further. Throughout, graphs are finite multigraphs. A 2-cycle means two parallel edges joining the same pair of vertices.

\begin{definition}
A matroid is called series--parallel if it is associated with a series--parallel graph.
\end{definition}

\begin{definition}\label{sp_def}
A graph is called series--parallel if it can be obtained from a loop or a bridge by a sequence of series and parallel extensions, i.e. subdivisions and duplications of edges.
\end{definition}

The series--parallel matroids referred to in \cite{A115594} are associated with the so-called \emph{2-connected series--parallel graphs} (further, simply \emph{2-sp-graphs}).

\begin{definition}\label{2sp_def}
A graph is called 2-sp-graph if it can be obtained from a 2-cycle by a sequence of series and parallel extensions.
\end{definition}

The rank of a graphical matroid equals the number of edges in every spanning tree of the associated graph. In the case of series--parallel graphs the rank is one less than the number of vertices.

This graph-theoretic context led us to the following construction.

Let $G$ be a 2-sp-graph with $k+1$ vertices and $n$ edges. Choose arbitrary directions for its edges. Let $B\in\R^{(k+1) \times n}$ be the incidence matrix, and let $W=\diag(w_e)\in\R^{n \times n}$ be a diagonal matrix of positive edge weights. We consider the scaled star space
\[
   U_{G,W}=\col\bigl(W^{1/2}B^T\bigr)\subset\R^n.
\]
The choice of edge directions only changes coordinate signs, and therefore does not affect principal angles to coordinate subspaces.

In the next section, for every 2-sp-graph $G$ we define a particular system of \emph{graph-induced edge weights} which make the scaled star space $U_{G,W}$ at least $\arccos(1/\sqrt n)$ away from every coordinate $k$-subspace:
\begin{equation}\label{geqacos}
   \theta_{\max}(U_{G,W},E_I) \geq \arccos \frac{1}{\sqrt n} \quad \text{for all } E_I.
\end{equation}
More specifically, with these edge weights the left-hand side of \eqref{geqacos} is either $\arccos(1/\sqrt n)$ or $\pi/2$, depending on whether the $k$-edge subgraph corresponding to $E_I$ is a spanning tree or not.

The final section proves that graph-induced edge weights are the unique positive weights, up to a common factor, satisfying \eqref{geqacos}.

\begin{remark}
By using matrix CS-decomposition (see \cite{Golub2013}), one can show that
\[
\theta_{\max}(U,E_I) = \theta_{\max}(U^\perp,E_I^\perp)
\]
for every $k$-dimensional subspace $U$. Therefore, the function \eqref{tf} takes the same values as the analogous function defined on $(n-k)$-dimensional subspaces of $\R^n$. In terms of 2-sp-graphs, this duality is realized by taking a graph dual (which is also a 2-sp-graph) with the inverted edge weights. The symmetry in the rows of Table \ref{tab} reflects this fact.
\end{remark}

\section{Graph-induced edge weights and coefficients}

For convenience, we adopt the following equivalent point of view on series--parallel graphs, which in particular justifies their second name -- \emph{two-terminal series--parallel graphs}.

\begin{definition}\label{tsp_def}
A two-terminal series--parallel graph with terminals $l$ and $r$ is either one edge $(l,r)$ or is obtained by one of the following operations:
\begin{enumerate}
\item a parallel composition
\[
   \Gamma=\Gamma_1\parallel\cdots\parallel \Gamma_p,
\]
where all left terminals are identified and all right terminals are identified;
\item a series composition
\[
   \Gamma=\Gamma_1\circ\cdots\circ \Gamma_s,
\]
where the right terminal of $\Gamma_i$ is identified with the left terminal of
$\Gamma_{i+1}$, $i = 1, \ldots, s-1$, and the same is true for the left terminals of $\Gamma$ and $\Gamma_1$ and the right terminals of $\Gamma$ and $\Gamma_s$.
\end{enumerate}
\end{definition}

Further, we will use the agreement that series and parallel compositions in the given recurrent definition alternate. For 2-sp-graphs we deal with, the final composition -- and thus the first decomposition -- is always parallel.

According to Lemma 9 from \cite{Eppstein1992}, if $(l, r)$ is an edge in a 2-sp-graph, $l$ and $r$ can be treated as the terminals.

\begin{figure}[H]
\centering
\includegraphics[width=12cm]{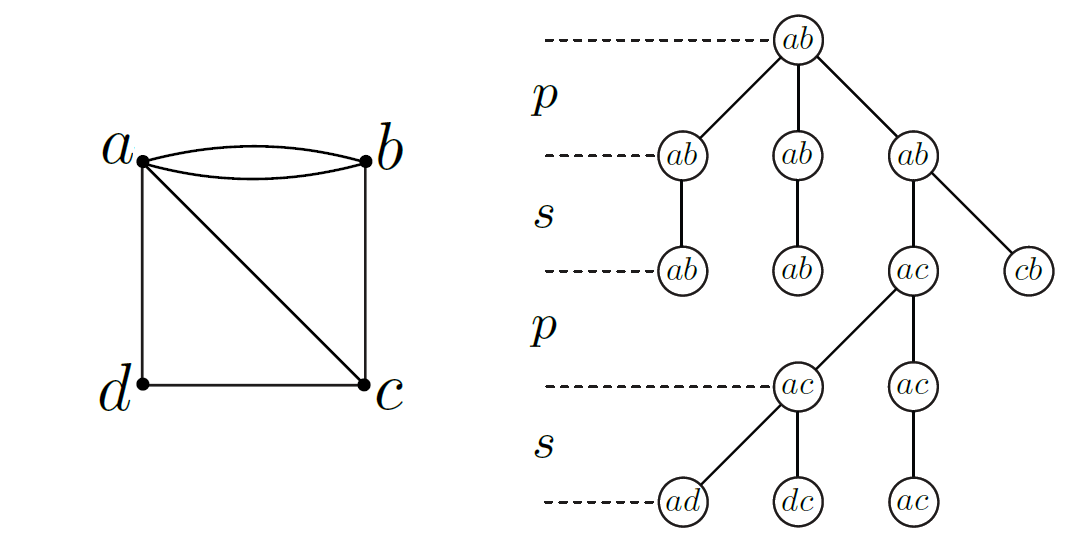}
\captionsetup{justification=centering}
\caption{ A 2-sp-graph and its admissible decomposition tree for terminal vertices $a$ and $b$. Nodes are labeled with their terminal pairs. }
\label{fig:quad}
\end{figure}

Definition \ref{tsp_def} introduces a \emph{decomposition tree} of a two-terminal series--parallel graph, whose root corresponds to the whole graph and all the other nodes correspond to the subgraphs given by series--parallel decomposition. According to the agreement above, when moving away from the root to a leaf, the types of decompositions alternate starting from the parallel. Also, if a chain of such decompositions ends with a parallel one, we will agree to complete it by a dummy series decomposition (see Figure \ref{fig:quad}). Decomposition trees satisfying these agreements will be called \emph{admissible}. A terminal choice determines an admissible decomposition tree uniquely, up to reduction of intermediate dummy decompositions with subsequent merging of neighboring same-type decompositions.

In what follows, $|\cdot|$ applied to a graph denotes its number of edges. In particular, $|G| = n$. Also, let
\[
   \varphi(x)=x(n-x).
\]

\begin{definition}
Given a 2-sp-graph $G$, fix a terminal pair and thus an admissible decomposition tree. For an edge $e$, let
\[
   G\supset \Gamma_1\supseteq \Gamma_2\supseteq\cdots\supseteq \Gamma_{2d}=\{e\}
\]
be the corresponding root-to-leaf chain. The graph-induced weight of $e$ is
\begin{equation}\label{eq:graph-induced-weight}
   w_e=\frac{\varphi(|\Gamma_1|)}{\varphi(|\Gamma_2|)}
        \frac{\varphi(|\Gamma_3|)}{\varphi(|\Gamma_4|)}
        \cdots
        \frac{\varphi(|\Gamma_{2d-1}|)}{\varphi(|\Gamma_{2d}|)}.
\end{equation}
\end{definition}

\begin{figure}[H]
\centering
\includegraphics[width=12cm]{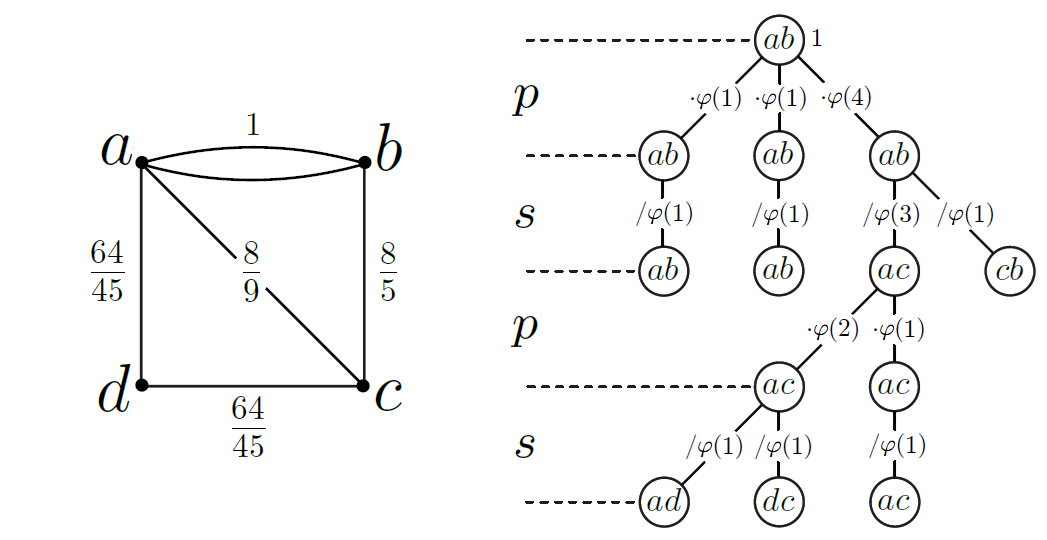}
\captionsetup{justification=centering}
\caption{ A 2-sp-graph and its graph-induced edge weights. }
\label{fig:quad_w}
\end{figure}

Figure \ref{fig:quad_w} illustrates formula \eqref{eq:graph-induced-weight}.

Dummy decompositions correspond to $\Gamma_i = \Gamma_{i+1}$, whose factors in \eqref{eq:graph-induced-weight} can be canceled. Therefore, graph-induced weights are determined exclusively by a 2-sp-graph and its terminals. Moreover, the following property holds.

\begin{proposition}
\label{prop:weight-terminal-invariance}
Changing the terminals changes all graph-induced weights by one common positive
factor.
\end{proposition}
\begin{proof}
It suffices to compare an arbitrary terminal choice with the terminal choice given by the endpoints of one fixed edge $a \in G$.

Denote the original decomposition tree by $\mathcal T$. Fix $e \in G$, and let $u$ be the lowest common ancestor of the leaves $a$ and $e$ in $\mathcal T$. In the quotient $w_e^{\mathcal T}/w_a^{\mathcal T}$, all factors from the root down to $u$ cancel. The factors from $u$ down to $e$ remain unchanged. Rewrite the factors from $u$ down to $a$ backwards reversing their exponents. Since $\varphi(|G\setminus\Gamma|) = \varphi(|\Gamma|)$, each node in that formula can be replaced by its complementary two-terminal subgraph. This yields formula \eqref{eq:graph-induced-weight} for a root-to-$e$ chain in an admissible decomposition tree $\mathcal T^a$ corresponding to the endpoints of $a$. Therefore,
\[
   \frac{w_e^{\mathcal T}}{w_a^{\mathcal T}} = w_e^{\mathcal T^a}, \quad e\in G.
\]

Finally, two arbitrary terminal choices can both be compared with the same terminal choice given by the endpoints of the edge $a$. This proves the proposition.
\end{proof}

\begin{remark}
Since every parallel composition of $2$-terminal 2-sp-graphs corresponds to a series composition of their duals and vice versa, the graph-induced edge weights of a 2-sp-graph are the reciprocals of the graph-induced edge weights of its dual, up to a common positive factor.
\end{remark}

We also define coefficients depending on a spanning tree.

\begin{definition}
Given a directed 2-sp-graph $G$ and its spanning tree $\tau$, fix a terminal pair with an admissible decomposition tree. For a node $\Gamma$ of this decomposition tree, set
\[
   \psi(\Gamma)=
   \begin{cases}
      n-|\Gamma|, & \tau\cap \Gamma\text{ connects the two terminals of } \Gamma,\\
      -|\Gamma|, & \text{otherwise}.
   \end{cases}
\]
For $e\in G$, with root-to-leaf chain
\[
   G\supset \Gamma_1\supseteq \Gamma_2\supseteq\cdots\supseteq \Gamma_{2d}=\{e\},
\]
define
\begin{equation}
\label{eq:induced-coefficient}
   y_e=\eps(e)
   \frac{\psi(\Gamma_1)}{\psi(\Gamma_2)}
   \frac{\psi(\Gamma_3)}{\psi(\Gamma_4)}
   \cdots
   \frac{\psi(\Gamma_{2d-1})}{\psi(\Gamma_{2d})},
\end{equation}
where $\eps(e)=1$ if the direction of $e$ agrees with the ordered terminals inherited at the leaf $\{e\}$, and $\eps(e)=-1$ otherwise. We call the coefficients $y_e$ $(G,\tau)$-induced.
\end{definition}

\begin{figure}[H]
\centering
\includegraphics[width=12cm]{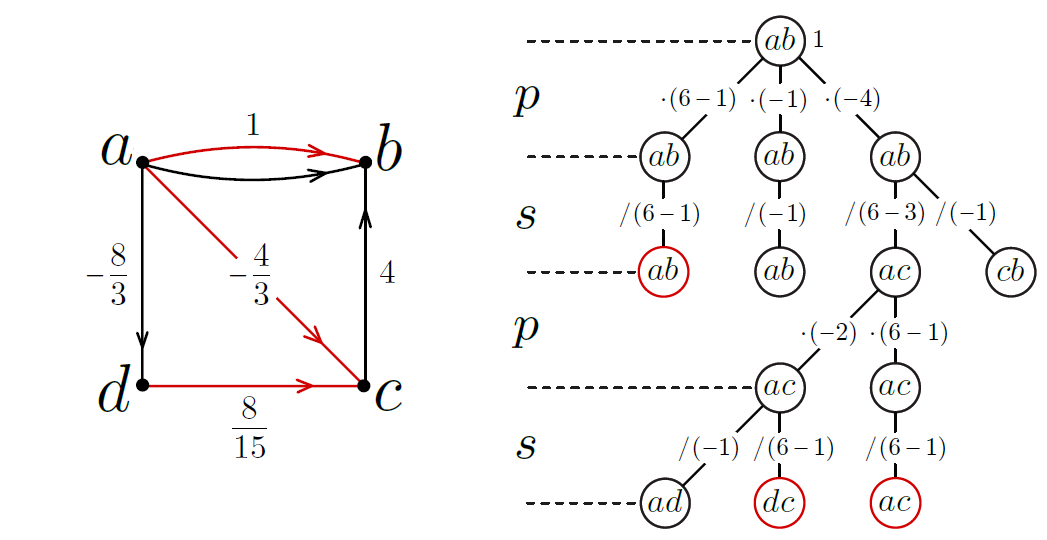}
\captionsetup{justification=centering}
\caption{ A directed 2-sp-graph, its spanning tree (in red), and the induced edge coefficients. }
\label{fig:quad_c}
\end{figure}

Figure \ref{fig:quad_c} illustrates formula \eqref{eq:induced-coefficient}.

Similar to graph-induced weights, coefficients \eqref{eq:induced-coefficient} are invariant under the choice of an admissible decomposition tree and determined exclusively by a 2-sp-graph and its terminals. Also, similarly, the following property holds.

\begin{proposition}
\label{prop:coeff-terminal-invariance}
Let $\mathcal T$ and $\mathcal S$ be two choices of admissible decomposition trees with the corresponding ordered terminals for the same directed 2-sp-graph $G$. Then there is a scalar $c \ne 0$ such that
\[
   y_e^{\mathcal S}=c\,y_e^{\mathcal T}, \quad e \in G.
\]
\end{proposition}

\begin{proof}
For a node $\Gamma$ of a decomposition tree $\mathcal T$, define
\[
   \chi_{\mathcal T}(\Gamma)=
   \begin{cases}
      1, & \tau\cap\Gamma\text{ connects the two terminals of }\Gamma,\\
      0, & \text{otherwise}.
   \end{cases}
\]
Then
\begin{equation}
\label{eq:psi-chi-combined}
   \psi_{\mathcal T}(\Gamma)=n\chi_{\mathcal T}(\Gamma)-|\Gamma|.
\end{equation}
If
\[
   G\supset \Gamma_1\supseteq\Gamma_2\supseteq\cdots\supseteq \Gamma_{2d}=\{e\}
\]
is the root-to-leaf chain, put $\alpha_i=(-1)^{i+1}$. Then
\[
   y_e^{\mathcal T}=\eps_{\mathcal T}(e)
      \prod_{i=1}^{2d}\psi_{\mathcal T}(\Gamma_i)^{\alpha_i}.
\]

It suffices to compare an arbitrary terminal choice with the terminal choice
given by the endpoints of one fixed edge $a \in G$. Choose those endpoints as terminals, ordered so that $a$ is positively oriented.

Similar to the proof of Proposition \ref{prop:weight-terminal-invariance}, fix $e \in G$, and let $u$ be the lowest
common ancestor of the leaves $a$ and $e$ in $\mathcal T$. In the quotient
$y_e^{\mathcal T}/y_a^{\mathcal T}$, all factors from the root down to $u$
cancel. The old branch from $u$ down to $e$ remains unchanged, while the old
branch from $u$ down to $a$ is read backwards; each node on that branch is
replaced by its complementary two-terminal subgraph, and its exponent is
reversed. This is a root-to-$e$ chain in an admissible decomposition tree $\mathcal T^a$ corresponding to the endpoints of $a$. We prove
\begin{equation}\label{eq:terminal-invariance-ratio}
   y_e^{\mathcal T^a}=\frac{y_e^{\mathcal T}}{y_a^{\mathcal T}}, \quad e\in G.
\end{equation}

Let $\Gamma$ be a node on the old $a$-branch, and let $\overline\Gamma$ be the
complementary node that appears after rerooting. Then
\[
   |\overline\Gamma|=n-|\Gamma|.
\]
Moreover,
\begin{equation}
\label{eq:chi-flip}
   \chi_{\mathcal T^a}(\overline\Gamma)=1-\chi_{\mathcal T}(\Gamma).
\end{equation}
Indeed, if $p,q$ are the two boundary terminals of $\Gamma$, the unique
$\tau$-path from $p$ to $q$ lies either inside $\Gamma$ or inside the
complementary two-terminal piece. Both alternatives cannot occur, since that
would create a cycle in $\tau$, and neither can fail, since $\tau$ is
connected. Hence the connecting status flips under complementation.

Using \eqref{eq:psi-chi-combined} and \eqref{eq:chi-flip}, we get
\begin{equation}
\label{eq:psi-flip}
   \psi_{\mathcal T^a}(\overline\Gamma)
   =n(1-\chi_{\mathcal T}(\Gamma))-(n-|\Gamma|)
   =-\psi_{\mathcal T}(\Gamma).
\end{equation}
Nodes on the old $e$-branch are not complemented, and their $\psi$-factors are
unchanged.

Let $N(a,e)$ be the number of complemented nodes on the new root-to-$e$ chain.
The product of $\psi$-factors in $y_e^{\mathcal T^a}$ is therefore
\[
   (-1)^{N(a,e)}
   \frac{
      \prod_{\Gamma\in\chain_{\mathcal T}(e)}
      \psi_{\mathcal T}(\Gamma)^{\alpha_{\mathcal T}(\Gamma)}
   }{
      \prod_{\Gamma\in\chain_{\mathcal T}(a)}
      \psi_{\mathcal T}(\Gamma)^{\alpha_{\mathcal T}(\Gamma)}
   }.
\]
The orientation sign changes by the same factor:
\[
   \eps_{\mathcal T^a}(e)=(-1)^{N(a,e)}
      \frac{\eps_{\mathcal T}(e)}{\eps_{\mathcal T}(a)}.
\]
Multiplying these two relations, the factors $(-1)^{N(a,e)}$ cancel, and
\eqref{eq:terminal-invariance-ratio} follows.

Finally, two arbitrary terminal choices can both be compared with the same terminal choice given by the endpoints of the edge $a$. This proves the proposition.
\end{proof}

\section{The main theorem}

\begin{theorem}\label{thm:main-construction}
Let $G$ be a directed 2-sp-graph with $k+1$ vertices and $n$ edges. Let $B\in\R^{(k+1) \times n}$ be the incidence matrix and $W\in\R^{n \times n}$ be the diagonal matrix of the $G$-induced edge weights. Then
\[
   U=\col(W^{1/2}B^T)\subset\R^n
\]
deviates from every coordinate $k$-subspace by either $\arccos(1/\sqrt n)$ or $\pi/2$.
\end{theorem}

For convenience, the proof is presented in three parts. The first, preparatory part introduces the main objects and the notation used later. The second part proves that subspace $U$ deviates from every coordinate $k$-subspace by at least $\arccos(1/\sqrt n)$. The third part strengthens this claim to that presented in the theorem.

\subsection*{Proof (1/3)}

We will use the two projection matrices onto subspace $U$ -- the orthogonal projection
\[
   P=W^{1/2}B^TL^+BW^{1/2},
\]
where $L^+$ is the Moore--Penrose inverse of the weighted graph Laplacian $L=BWB^T$,
and the diagonally similar matrix
\[
   Y=W^{1/2}PW^{-1/2}=WB^TL^+B,
\]
also known as the transfer-current matrix of the directed weighted graph $G$ (see \cite{Kassel2015}). Since the similarity is diagonal, corresponding principal blocks of $P$ and $Y$ are also similar and thus have the same eigenvalues.

By Theorem 1.16 from \cite{Grimmett2018}, $Y$ has the following explicit form
\begin{equation}\label{Ygen}
Y = \frac{1}{T} \, \begin{blockarray}{cccccc}
  & e &   & f &   \\
\begin{block}{[ccccc]c}
    & \vdots  &   & \vdots &   &   \\
\,\, \cdots & T_{e+} & \cdots & T^{\rightrightarrows}_{e+f-}-T^{\rightleftarrows}_{e+f-} & \cdots \,\, & e \\
    & \vdots &   & \vdots &   &   \\
\,\, \cdots & T^{\rightrightarrows}_{f+e-}-T^{\rightleftarrows}_{f+e-} & \cdots & T_{f+} & \cdots \,\, & f \\
    & \vdots &   & \vdots &   &   \\
\end{block}
\end{blockarray},
\end{equation}
where $T$ denotes the following weighted sum over spanning trees of $G$
\[
   T=\sum_t\prod_{e\in t}w_e,
\]
$T_{e+}$ denotes the similar weighted sum over spanning trees containing edge $e$, and $T^{\rightrightarrows}_{e+f-}$ and $T^{\rightleftarrows}_{e+f-}$ are the similar weighted sums over spanning trees containing $e$ for which
the unique path between the ordered endpoints of $f$ passes through $e$ in the positive and negative direction, respectively.

By Theorem 3 from \cite{Duffin1965}, any two edges of a 2-sp-graph oriented coherently in a cycle cannot be presented in a different cycle in the opposite way (and vice versa). This yields the following simplified version of \eqref{Ygen}
\begin{equation}\label{Ysp}
Y = \frac{1}{T} \, \begin{blockarray}{cccccc}
  & e &   & f &   \\
\begin{block}{[ccccc]c}
    & \vdots  &   & \vdots &   &   \\
\,\, \cdots & T_{e+} & \cdots & \pm T_{e+f-} & \cdots \,\, & e \\
    & \vdots &   & \vdots &   &   \\
\,\, \cdots & \pm T_{f+e-} & \cdots & T_{f+} & \cdots \,\, & f \\
    & \vdots &   & \vdots &   &   \\
\end{block}
\end{blockarray},
\end{equation}
where $T_{e+f-}$ has the same meaning as given above but without taking into account edge directions. The latter affects the sign before the expression.

We use analogous notation for weighted sums over spanning trees of a subgraph $\Gamma \subset G$: $T(\Gamma)$. For a single edge $e$, $T(\{e\}) = w_e$.

Also, for a two-terminal graph $\Gamma$, let
\[
   F(\Gamma)=\sum_f\prod_{e\in f}w_e
\]
be the weighted sum over two-component forests rooted at the two terminals.
For a single edge $e$, $F(\{e\})=1$.

By Theorem 1 from \cite{Jan1992}, the following product rules hold
\begin{equation}\label{eq:serial-product-rules}
\begin{split}
   T(\Gamma_1\circ\cdots\circ\Gamma_s) &=\prod_{j=1}^sT(\Gamma_j), \\
   \frac{F(\Gamma_1\circ\cdots\circ\Gamma_s)}{T(\Gamma_1\circ\cdots\circ\Gamma_s)}
      &=\sum_{j=1}^s\frac{F(\Gamma_j)}{T(\Gamma_j)},
\end{split}
\end{equation}
\begin{equation}\label{eq:parallel-product-rules}
\begin{split}
   F(\Gamma_1\parallel\cdots\parallel\Gamma_p) &=\prod_{i=1}^p F(\Gamma_i), \\
   \frac{T(\Gamma_1\parallel\cdots\parallel\Gamma_p)}
        {F(\Gamma_1\parallel\cdots\parallel\Gamma_p)}
      &=\sum_{i=1}^p\frac{T(\Gamma_i)}{F(\Gamma_i)}.
\end{split}
\end{equation}

\subsection*{Proof (2/3)}

We now prove that every diagonal $k$-block of orthogonal projection $P$ is either singular or has an eigenvalue equal to $1/n$.

If $I\subset E(G)$ has $|I|=k$ but is not a spanning tree, then the
incidence block $B_I$ has rank less than $k$. Thus $P_I$ is singular, and the corresponding largest principal angle is $\pi/2$.

It remains to consider the case where $I = \tau$ is a spanning tree. We prove
\begin{equation}
\label{eq:Ytau-eigen-main}
   Y_\tau y_\tau=\frac1n y_\tau,
\end{equation}
where $y_\tau = (y_e)_{e \in I}$ is the restriction of the vector of $(G,\tau)$-induced coefficients to the edges of $\tau$. This implies that each spanning-tree block has an eigenvalue equal to $1/n$.

Fix an arbitrary edge $e_0\in\tau$. It is enough to prove the $e_0$-coordinate
of \eqref{eq:Ytau-eigen-main}. By Propositions
\ref{prop:weight-terminal-invariance} and \ref{prop:coeff-terminal-invariance},
we may choose the endpoints of $e_0$ as terminals without loss of generality:
the weights are only multiplied by one common positive factor, and
$y_\tau$ is only multiplied by one common nonzero factor. We order the
terminals as $(l,r)$ so that $e_0$ is positively oriented. By changing
coordinate signs, we may assume that each edge direction agrees with its
inherited left-to-right terminal orientation. In this normalization,
\[
   w_{e_0}=1,
   \quad
   y_\tau(e_0)=1.
\]

With these terminals the first parallel decomposition has the form
\[
   G=\{e_0\}\parallel H, \quad H = H_1 \parallel \cdots \parallel H_p,
\]
where $H$ is the complementary two-terminal series--parallel graph. Since
$e_0\in\tau$, the restriction $\tau\cap H$ is a two-component forest separating
$l$ and $r$.

This yields
\[
   T_{0+}=F(H),
   \quad
   T_{0-}=T(H),
   \quad
   T=T_{0+}+T_{0-},
\]
where $T_{0+}$ and $T_{0-}$ are the sums over spanning trees containing and not containing $e_0$, respectively.

The transfer-current formula \eqref{Ysp} gives
\begin{equation}
\label{eq:transfer-row-main}
   \bigl(TY_\tau y_\tau\bigr)(e_0)
   =T_{0+}+
     \sum_{e\in\tau\cap H}T_{0+e-}\,y_\tau(e).
\end{equation}
Here $T_{0+e-}$ denotes the weighted sum of spanning trees containing $e_0$ for
which the unique path between the endpoints of $e$ passes through $e_0$.
Thus the desired coordinate identity follows once we prove
\begin{equation}
\label{eq:root-goal-main}
   \sum_{e\in\tau\cap H}T_{0+e-}\,y_\tau(e)
   =\frac{T_{0-}-(n-1)T_{0+}}{n}.
\end{equation}
Indeed, substituting \eqref{eq:root-goal-main} into
\eqref{eq:transfer-row-main} gives
\[
   \bigl(TY_\tau y_\tau\bigr)(e_0)
   =T_{0+}+\frac{T_{0-}-(n-1)T_{0+}}{n}
   =\frac{T}{n}
   =\frac{T}{n}y_\tau(e_0).
\]

Consider an admissible decomposition tree of $H$. We call the root $H$, and
every node reached after a serial decomposition, a \emph{parallel node}. The remaining nodes will be called \emph{serial}. Thus a
non-leaf parallel node $A$ has a parallel decomposition
\[
   A=S_1\parallel\cdots\parallel S_p,
\]
where the $S_i$ are serial nodes. A serial node $S$ has a serial decomposition
\[
   S=A_1\circ\cdots\circ A_s,
\]
where the $A_j$ are parallel nodes. Leaves are single-edge parallel nodes.

For a parallel node $A$, define the accumulated weight factor $w(A)$ by
\[
   w(H)=1,
\]
and, whenever
\[
   A=S_1\parallel\cdots\parallel S_p,
   \quad
   S_i=A_{i1}\circ\cdots\circ A_{is_i},
\]
by
\begin{equation}
\label{eq:w-main}
   w(A_{ij})=
   w(A) \, \frac{|S_i|(n-|S_i|)}{|A_{ij}|(n-|A_{ij}|)}.
\end{equation}
For a leaf $\{e\}$, this gives precisely the graph-induced edge weight:
\begin{equation}
\label{eq:w-leaf-main}
   w(\{e\})=w_e.
\end{equation}

For a parallel node $A$, define $T_{0+A-}$ to be the weighted sum of spanning
trees of $G$ containing $e_0$ such that the unique path between the two
terminals of $A$ passes through $e_0$. Equivalently, inside $A$ one uses a
separating two-component forest, and the rest of $G$ is completed compatibly.
In particular,
\begin{equation}
\label{eq:T-root-main}
   T_{0+H-} = F(H) = T_{0+}.
\end{equation}
If $A_{ij}$ is a serial child of the parallel child $S_i$ of $A$, then
\begin{equation}
\label{eq:T-main}
   T_{0+A_{ij}-}=
   T_{0+A-}\,
   \frac{T(S_i)}{F(S_i)}
   \frac{F(A_{ij})}{T(A_{ij})}.
\end{equation}
Indeed, in the event counted by $T_{0+A_{ij}-}$, the branch $S_i$ is forced to
have $A_{ij}$ separated and all other serial children connected. Dividing this
specified contribution by the former contribution $F(S_i)$ gives
\eqref{eq:T-main}.

The fixed tree $\tau$ gives each two-terminal subgraph one of two boundary
states: connected, if its two terminals are connected inside the subgraph, and
separated otherwise. At a parallel node, a connected state means that exactly
one parallel child is connected and all other parallel children are separated;
a separated state means that all parallel children are separated. At a serial
node, a connected state means all serial children are connected; a separated
state means exactly one serial child is separated and all others are connected.
Whenever such a child is distinguished, we label it first.

\begin{figure}[H]
\centering
\includegraphics[width=13cm]{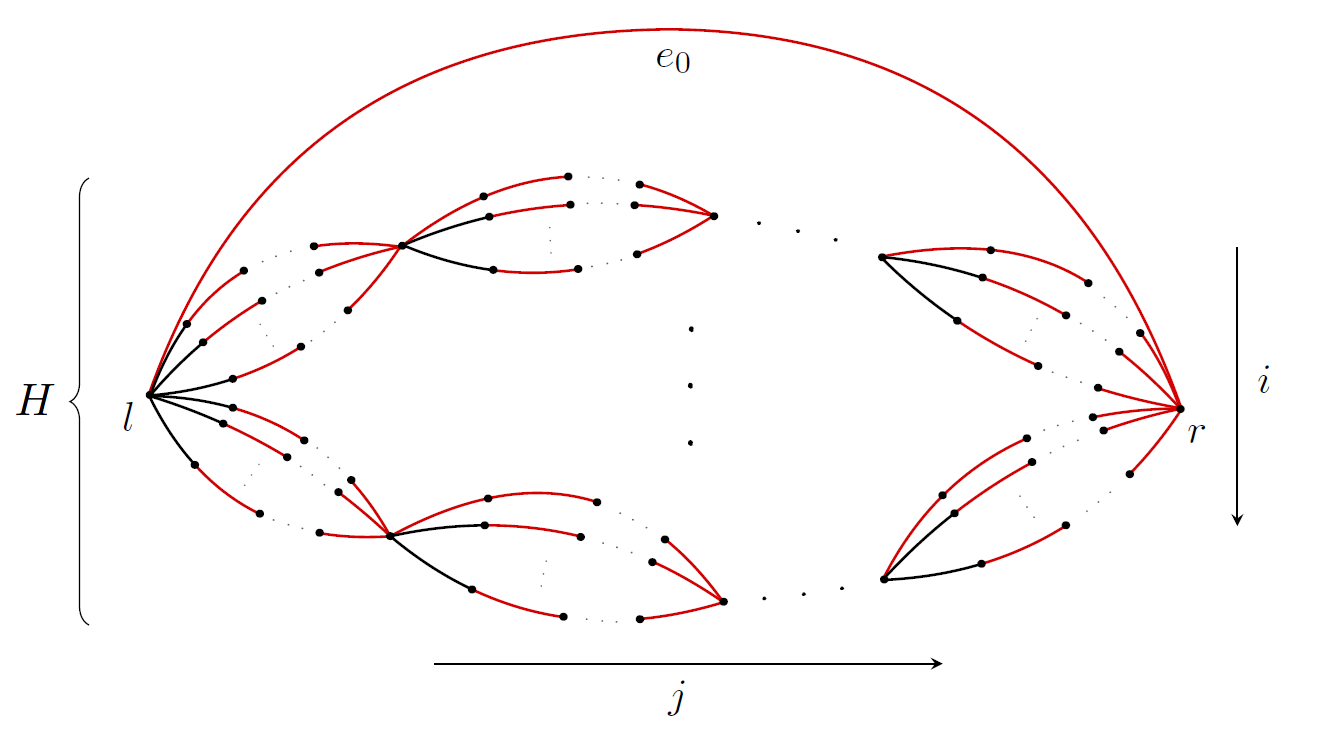}
\captionsetup{justification=centering}
\caption{ Graph $G$ and the layout of $\tau$ (in red). }
\label{fig:SPtree_upd}
\end{figure}

Figure \ref{fig:SPtree_upd} illustrates this agreement.

\vspace*{-0.25cm}

For a parallel node $A$, define
\[
   C^{\mathrm{sep}}(A), \quad C^{\mathrm{con}}(A)
\]
to be the total contribution of the edges of $\tau\cap A$ to the sum in
\eqref{eq:root-goal-main}, with all coefficient factors above $A$ omitted,
under the separated and connected boundary states, respectively. For a serial
node $S$, define analogous quantities
\[
   D^{\mathrm{sep}}(S),\quad D^{\mathrm{con}}(S).
\]
For a leaf $A=\{e\}$,
\begin{equation}
\label{eq:C-base-main}
   C^{\mathrm{sep}}(\{e\})=0,
   \quad
   C^{\mathrm{con}}(\{e\})=T_{0+e-}.
\end{equation}
For a parallel node $A=S_1\parallel\cdots\parallel S_p$,
\begin{align}
\label{eq:Csep-rec-main}
   C^{\mathrm{sep}}(A)
      &=\sum_{i=1}^p(-|S_i|)D^{\mathrm{sep}}(S_i),\\
\label{eq:Ccon-rec-main}
   C^{\mathrm{con}}(A)
      &=(n-|S_1|)D^{\mathrm{con}}(S_1)
        +\sum_{i=2}^p(-|S_i|)D^{\mathrm{sep}}(S_i),
\end{align}
where $S_1$ is the connected child in the connected case. For a serial node
$S=A_1\circ\cdots\circ A_s$,
\begin{align}
\label{eq:Dcon-rec-main}
   D^{\mathrm{con}}(S)
      &=\sum_{j=1}^s\frac1{n-|A_j|}C^{\mathrm{con}}(A_j),\\
\label{eq:Dsep-rec-main}
   D^{\mathrm{sep}}(S)
      &=\frac1{-|A_1|}C^{\mathrm{sep}}(A_1)
        +\sum_{j=2}^s\frac1{n-|A_j|}C^{\mathrm{con}}(A_j),
\end{align}
where $A_1$ is the separated child in the separated case. These are exactly
the factors in \eqref{eq:induced-coefficient}: parallel levels multiply
by $\psi(\Gamma)$, and serial levels divide by $\psi(\Gamma)$. Hence
\begin{equation}
\label{eq:root-sum-C-main}
   \sum_{e\in\tau\cap H}T_{0+e-}\,y_\tau(e)=C^{\mathrm{sep}}(H),
\end{equation}
because $H$ is separated when $e_0\in\tau$.

We now evaluate these quantities in closed form.

\begin{lemma}\label{lem:closed-main}
For every parallel node $A$,
\begin{align}\label{eq:Csep-closed-main}
   C^{\mathrm{sep}}(A)
      &=\frac{T_{0+A-}}{n}
        \left(
        \frac{T(A)}{w(A)F(A)}-|A|
        \right),\\
\label{eq:Ccon-closed-main}
   C^{\mathrm{con}}(A)
      &=\frac{T_{0+A-}}{n}
        \left(
        \frac{T(A)}{w(A)F(A)}+n-|A|
        \right).
\end{align}
\end{lemma}

\begin{proof}
We argue by induction upward from the leaves. If $A=\{e\}$, then
$T(A)=w_e$, $F(A)=1$, and $w(A)=w_e$. Hence $T(A)/(w(A)F(A))=1$ and $|A|=1$, so the right-hand sides of \eqref{eq:Csep-closed-main} and \eqref{eq:Ccon-closed-main} are $0$ and $T_{0+e-}$, in agreement with \eqref{eq:C-base-main}.

Now let
\[
   A=S_1\parallel\cdots\parallel S_p,
   \quad
   S_i=A_{i1}\circ\cdots\circ A_{is_i},
\]
and assume the formulas are known for all parallel children $A_{ij}$.
First compute $D^{\mathrm{con}}(S_i)$. From \eqref{eq:Dcon-rec-main} and the
induction hypothesis,
\[
   D^{\mathrm{con}}(S_i)
   =\sum_j\frac1{n-|A_{ij}|}
      \frac{T_{0+A_{ij}-}}{n}
      \left(
      \frac{T(A_{ij})}{w(A_{ij})F(A_{ij})}
      +n-|A_{ij}|
      \right).
\]
Using \eqref{eq:w-main} and \eqref{eq:T-main}, the part involving
$w(A_{ij})$ contributes
\[
   \frac{T_{0+A-}}{n}
   \frac{T(S_i)}{F(S_i)}
   \sum_j\frac1{w(A_{ij})(n-|A_{ij}|)}
   =
   \frac{T_{0+A-}}{n}
   \frac{T(S_i)}
        {w(A)F(S_i)(n-|S_i|)}.
\]
Indeed,
\[
   \frac1{w(A_{ij})(n-|A_{ij}|)}
   =\frac{|A_{ij}|}{w(A)|S_i|(n-|S_i|)},
\]
and $\sum_j|A_{ij}|=|S_i|$. The remaining part contributes
\[
   \frac{T_{0+A-}}{n}
   \frac{T(S_i)}{F(S_i)}
   \sum_j\frac{F(A_{ij})}{T(A_{ij})}
   =\frac{T_{0+A-}}{n}
\]
by \eqref{eq:serial-product-rules}. Therefore
\begin{equation}
\label{eq:Dcon-closed-main}
   D^{\mathrm{con}}(S_i)
   =\frac{T_{0+A-}}{n}
      \left(
      1+
      \frac{T(S_i)}
           {w(A)F(S_i)(n-|S_i|)}
      \right).
\end{equation}

Next compute $D^{\mathrm{sep}}(S_i)$. Label the separated serial child by
$A_{i1}$. From \eqref{eq:Dsep-rec-main},
\begin{align*}
   D^{\mathrm{sep}}(S_i)
   &=\frac1{-|A_{i1}|}
      \frac{T_{0+A_{i1}-}}{n}
      \left(
      \frac{T(A_{i1})}{w(A_{i1})F(A_{i1})}
      -|A_{i1}|
      \right) \\
   &\quad+
      \sum_{j=2}^{s_i}\frac1{n-|A_{ij}|}
      \frac{T_{0+A_{ij}-}}{n}
      \left(
      \frac{T(A_{ij})}{w(A_{ij})F(A_{ij})}
      +n-|A_{ij}|
      \right).
\end{align*}
The terms not involving $w(\cdot)$ again combine to $T_{0+A-}/n$. The terms
involving $w(\cdot)$ combine to
\[
   \frac{T_{0+A-}}{n}
   \frac{T(S_i)}{F(S_i)}
   \left(
      -\frac1{w(A_{i1})|A_{i1}|}
      +\sum_{j=2}^{s_i}\frac1{w(A_{ij})(n-|A_{ij}|)}
   \right).
\]
By \eqref{eq:w-main},
\[
   \frac1{w(A_{i1})|A_{i1}|}
   =\frac{n-|A_{i1}|}{w(A)|S_i|(n-|S_i|)},
\]
whereas for $j\ge2$,
\[
   \frac1{w(A_{ij})(n-|A_{ij}|)}
   =\frac{|A_{ij}|}{w(A)|S_i|(n-|S_i|)}.
\]
Since $\sum_{j=2}^{s_i}|A_{ij}|=|S_i|-|A_{i1}|$, the expression in parentheses
is $-1/(w(A)|S_i|)$. Thus
\begin{equation}
\label{eq:Dsep-closed-main}
   D^{\mathrm{sep}}(S_i)
   =\frac{T_{0+A-}}{n}
      \left(
      1-
      \frac{T(S_i)}
           {w(A)F(S_i)|S_i|}
      \right).
\end{equation}

Substitute \eqref{eq:Dsep-closed-main} into \eqref{eq:Csep-rec-main}. We get
\begin{align*}
   C^{\mathrm{sep}}(A)
   &=\sum_{i=1}^p(-|S_i|)D^{\mathrm{sep}}(S_i) \\
   &=\frac{T_{0+A-}}{n}
      \sum_{i=1}^p
      \left(
      -|S_i|+
      \frac{T(S_i)}{w(A)F(S_i)}
      \right) \\
   &=\frac{T_{0+A-}}{n}
      \left(
      \frac1{w(A)}
      \sum_{i=1}^p\frac{T(S_i)}{F(S_i)}-|A|
      \right) \\
   &=\frac{T_{0+A-}}{n}
      \left(
      \frac{T(A)}{w(A)F(A)}-|A|
      \right),
\end{align*}
where the last equality is \eqref{eq:parallel-product-rules}.

For the connected state, label the connected parallel child by $S_1$. Using
\eqref{eq:Dcon-closed-main}, \eqref{eq:Dsep-closed-main}, and
\eqref{eq:Ccon-rec-main},
\begin{align*}
   C^{\mathrm{con}}(A)
   &=(n-|S_1|)D^{\mathrm{con}}(S_1)
     +\sum_{i=2}^p(-|S_i|)D^{\mathrm{sep}}(S_i) \\
   &=\frac{T_{0+A-}}{n}
      \left(
      n-|S_1|+
      \frac{T(S_1)}{w(A)F(S_1)}
      \right) \\
   &\quad+
      \frac{T_{0+A-}}{n}
      \sum_{i=2}^p
      \left(
      -|S_i|+
      \frac{T(S_i)}{w(A)F(S_i)}
      \right) \\
   &=\frac{T_{0+A-}}{n}
      \left(
      \frac1{w(A)}
      \sum_{i=1}^p\frac{T(S_i)}{F(S_i)}+n-|A|
      \right) \\
   &=\frac{T_{0+A-}}{n}
      \left(
      \frac{T(A)}{w(A)F(A)}+n-|A|
      \right).
\end{align*}
This completes the induction.
\end{proof}

Apply Lemma \ref{lem:closed-main} to the root $H$ in the separated state. Since
\[
   T_{0+H-}=F(H)=T_{0+},
   \quad
   w(H)=1,
   \quad
   |H|=n-1,
   \quad
   T(H)=T_{0-},
\]
we obtain
\[
   C^{\mathrm{sep}}(H)
   =\frac{T_{0+}}{n}
      \left(\frac{T_{0-}}{T_{0+}}-(n-1)\right)
   =\frac{T_{0-}-(n-1)T_{0+}}{n}.
\]
Together with \eqref{eq:root-sum-C-main}, this proves
\eqref{eq:root-goal-main}, and hence the $e_0$-coordinate of
\eqref{eq:Ytau-eigen-main}. Since $e_0\in\tau$ was arbitrary,
\[
   Y_\tau y_\tau=\frac1n y_\tau.
\]

Thus every spanning-tree block $P_\tau$ has $1/n$ as an eigenvalue, while every
non-spanning-tree $k$-block is singular. Let $Q$ be an $n\times k$ matrix with
orthonormal columns spanning $U$. Then $P=QQ^T$, and the squared cosines of
the principal angles between $U$ and the coordinate subspace $E_I$ are the
eigenvalues of $Q_IQ_I^T=P_I$. Therefore, for every $I$ with $|I|=k$,
\[
   \cos^2\theta_{\max}(U,E_I)=\lambda_{\min}(P_I)\le\frac1n,
\]
which is equivalent to
\[
   \theta_{\max}(U,E_I)\ge\arccos\frac1{\sqrt n}.
\]

\subsection*{Proof (3/3)}

We now show that $1/n$ is the least eigenvalue of $P_{\tau}$, completing the proof of the theorem.

Delete one row of the incidence matrix $B$ and denote the resulting reduced
incidence matrix by $\widetilde B$. Put
\[
   A=W^{1/2}\widetilde B^{\,T}\in\R^{n\times k}.
\]
Then
\[
   P=A(A^TA)^{-1}A^T.
\]
Since $\tau$ is a spanning tree, the square block $A_\tau$ is invertible. Let
$\bar\tau=E(G)\setminus\tau$, and set
\[
   R=A_{\bar\tau}A_\tau^{-1}.
\]
After ordering the rows of $A$ as $\tau,\bar\tau$, we have
\[
   A^TA=A_\tau^T(I+R^TR)A_\tau.
\]
Consequently
\begin{equation}
\label{eq:K-inverse}
   K:=P_\tau=(I+R^TR)^{-1},
   \quad
   K^{-1}=I+R^TR.
\end{equation}

For a non-tree edge $f\in\bar\tau$, let $\operatorname{path}_\tau(f)$ be the
unique path in $\tau$ joining the endpoints of $f$. The row of $R$ indexed by
$f$ is the weighted fundamental-cycle row:
\begin{equation}
\label{eq:R-fundamental-row}
   R_{fe}=\begin{cases}
      \eta_{fe}\sqrt{w_f/w_e}, & e\in\operatorname{path}_\tau(f),\\
      0, & e\notin\operatorname{path}_\tau(f),
   \end{cases}
\end{equation}
where $\eta_{fe}\in\{\pm1\}$ is the incidence sign of $e$ in the fundamental
cycle $\operatorname{path}_\tau(f)\cup\{f\}$.

We need a sign lemma.

\begin{lemma}\label{lem:fundamental-cycle-sign}
For each $f\in\bar\tau$, the nonzero numbers
\[
   R_{fe}y_e,
   \quad e\in\operatorname{path}_\tau(f),
\]
all have the same sign.
\end{lemma}

\begin{proof}
Choose the endpoints of $f$ as the terminals, ordered so that $f$ is positively
oriented. With this terminal choice the first decomposition has the form
\[
   G=\{f\}\parallel \widetilde{H}.
\]
For every edge $e$ on $\operatorname{path}_\tau(f)$, every subgraph $\Gamma$
on the root-to-$e$ chain is traversed by the unique $\tau$-path between the two
terminals of $\Gamma$. Hence $\tau\cap\Gamma$ connects the terminals of
$\Gamma$, and every corresponding factor is
\[
   \psi(\Gamma)=n-|\Gamma|>0.
\]
Thus, for this terminal choice, the sign of the induced coefficient is exactly
the fundamental-cycle incidence sign, up to one global sign depending only on
$f$. By Proposition \ref{prop:coeff-terminal-invariance}, changing the
terminals rescales all induced coefficients by one nonzero scalar. Therefore
the same sign relation holds for the original coefficient vector. Since all
weights are positive, \eqref{eq:R-fundamental-row} gives the claim.
\end{proof}

Let $y=y_\tau$ be the induced coefficient vector and put
\[
   x=W_\tau^{-1/2}y.
\]
Because $Y_\tau=W_\tau^{1/2}P_\tau W_\tau^{-1/2}$ and
$Y_\tau y=(1/n)y$, we have
\begin{equation}
\label{eq:K-eigenvector}
   Kx=\frac1n x,
   \quad
   K^{-1}x=nx.
\end{equation}
The vector $x$ has the same sign pattern as $y$, and all its entries are
nonzero. Let
\[
   D=\diag(\sgn x_e)_{e\in\tau}.
\]
By Lemma \ref{lem:fundamental-cycle-sign}, every row of $RD$ is either
entrywise nonnegative or entrywise nonpositive. Therefore
\[
   DR^TRD=(RD)^T(RD)
\]
is entrywise nonnegative. From \eqref{eq:K-inverse},
\[
   M:=DK^{-1}D=I+DR^TRD
\]
is an entrywise nonnegative matrix. Conjugating \eqref{eq:K-eigenvector} by
$D$ gives
\[
   M|x|=n|x|,
\]
where $|x|$ is strictly positive. By the Perron--Frobenius theorem, an
entrywise nonnegative matrix with a strictly positive eigenvector of eigenvalue
$n$ has spectral radius $n$. Thus $\rho(M)=n$.

The matrices $M$ and $K^{-1}$ are similar. Since $K^{-1}$ is symmetric
positive definite, its spectral radius is its largest eigenvalue. Hence
\[
   \lambda_{\max}(K^{-1})=n,
\]
and therefore
\[
   \lambda_{\min}(K)=\frac1n.
\]
Finally, $Y_\tau$ is diagonally similar to $K=P_\tau$, so
$\lambda_{\min}(Y_\tau)=1/n$ as well.\qed

\section{Fixed-graph optimality and uniqueness of the weights}

We will prove that for a fixed $2$-sp-graph, the graph-induced edge weights are the unique positive weights, up to a common factor, that make the bound $1/n$ sharp.

Throughout this section $G=(V,E)$ is a fixed $2$-sp-graph with
\[
   |E|=n, \quad |V|=k+1,
\]
and $W=\diag(w_e)$ is an arbitrary positive edge weighting. Define
\[
   \alpha_G(W)= \max_{\substack{I\subseteq E\\ |I|=k}}\lambda_{\min}P_I(W).
\]

\subsection{Principal blocks and energy ratios}

For an edge set $S\subseteq E$ and a potential $u:V\to\R$, set
\[
   \Energy_S(u)=\sum_{xy\in S}w_{xy}\bigl(u(x)-u(y)\bigr)^2.
\]
If $\tau$ is a spanning tree, define
\[
   \rho_\tau(W)=
   \sup_{u\not\equiv\mathrm{const}}
   \frac{\Energy_{E\setminus\tau}(u)}{\Energy_\tau(u)}.
\]

\begin{lemma}\label{lem:block-energy-combined}
For every spanning tree $\tau$,
\[
   \lambda_{\min}P_\tau(W)=\frac1{1+\rho_\tau(W)}.
\]
Consequently,
\[
   \alpha_G(W)=\frac1{1+\min_{\tau\in\SpTrees(G)}\rho_\tau(W)},
\]
where $\SpTrees(G)$ denotes the set of spanning trees of $G$.
\end{lemma}

\begin{proof}
With $A=W^{1/2}\widetilde B^{\,T}$, as in the proof of Theorem
\ref{thm:main-construction}, we have $P=A(A^TA)^{-1}A^T$. Since $\tau$ is a
spanning tree, $A_\tau$ is invertible. With
\[
   R_\tau=A_{E\setminus\tau}A_\tau^{-1},
\]
we have
\[
   P_\tau(W)=(I+R_\tau^TR_\tau)^{-1}.
\]
Moreover,
\[
   \|A_\tau u\|^2=\Energy_\tau(u),
   \quad
   \|A_{E\setminus\tau}u\|^2=\Energy_{E\setminus\tau}(u).
\]
As $A_\tau$ is invertible, $\rho_\tau(W)=\|R_\tau\|^2$. Therefore
\[
   \lambda_{\min}P_\tau(W)
   =\frac1{1+\lambda_{\max}(R_\tau^TR_\tau)}
   =\frac1{1+\rho_\tau(W)}.
\]
If $I$ is not a spanning tree, then $P_I(W)$ is singular. Hence the maximum
over all $k$-element sets is the maximum over spanning trees, and the second
formula follows because $x\mapsto(1+x)^{-1}$ is decreasing.
\end{proof}

Thus $\alpha_G(W)=1/n$ is equivalent to
\[
   \min_{\tau\in\SpTrees(G)}\rho_\tau(W)=n-1.
\]

\subsection{A tree--forest certificate}

Fix an admissible series--parallel decomposition of $G$. For a two-terminal subgraph $H$ in the decomposition tree, write its terminals as $l_H,r_H$ and put
\[
   |H|=|E(H)|,
   \quad
   \Delta_H(u)=u(l_H)-u(r_H).
\]
A \emph{tree of $H$} means a spanning tree of $H$. A \emph{separating forest
of $H$} means a spanning forest with exactly two components, one containing
$l_H$ and one containing $r_H$.

For $S\subseteq E(H)$ define
\begin{equation}
\label{eq:Q-def-combined}
   Q_H(S;u)=\Energy_{E(H)\setminus S}(u)-(n-1)\Energy_S(u).
\end{equation}
This quantity uses the total edge count $n$ of the ambient graph.

\begin{lemma}\label{lem:harmonic-cauchy-combined}
Let $A_1,\ldots,A_s>0$ and
\[
   A=\left(\sum_{i=1}^s\frac1{A_i}\right)^{-1}.
\]
Then for any $x_1, \ldots, x_s \in \R$
\[
   \sum_{i=1}^sA_i x_i^2\ge
   A\left(\sum_{i=1}^s x_i\right)^2.
\]
\end{lemma}

\begin{proof}
This is Cauchy's inequality applied to
$\sum_i x_i=\sum_i(\sqrt{A_i}x_i)(1/\sqrt{A_i})$.
\end{proof}

\begin{lemma}\label{lem:one-positive-combined}
Let $B>0$ and $A_i>0$ for $i\ne j$. Suppose
\[
   D^{-1}=\frac1B-\sum_{i\ne j}\frac1{A_i}>0.
\]
Then
\[
   Bx_j^2-\sum_{i\ne j}A_i x_i^2
   \le D\left(\sum_i x_i\right)^2.
\]
\end{lemma}

\begin{proof}
Fix $S=\sum_i x_i$. If $x_j$ is fixed, Lemma
\ref{lem:harmonic-cauchy-combined} shows that the expression is maximized over
the remaining variables by minimizing $\sum_{i\ne j}A_i x_i^2$. With
\[
   A_*=\biggl(\sum_{i\ne j}\frac1{A_i}\biggr)^{-1},
\]
it is enough to maximize
\[
   Bx_j^2-A_*(S-x_j)^2.
\]
Since $D^{-1}=B^{-1}-A_*^{-1}>0$, this quadratic polynomial has finite maximum
\[
   \frac{BA_*}{A_*-B}S^2
   =\left(\frac1B-\frac1{A_*}\right)^{-1}S^2
   =DS^2.
\]
\end{proof}

\begin{proposition}\label{prop:tree-forest-combined}
Let $H$ be any proper two-terminal subgraph in the chosen decomposition of
$G$. Then there exists a number $c_H>0$ such that both
statements below hold.
\begin{enumerate}
\item There is a tree $T_H$ of $H$ satisfying
\begin{equation}\label{eq:tree-invariant-combined}
   Q_H(T_H;u)\le -(n-|H|)c_H\Delta_H(u)^2
   \quad\text{for all }u.
\end{equation}
\item There is a separating forest $F_H$ of $H$ satisfying
\begin{equation}\label{eq:forest-invariant-combined}
   Q_H(F_H;u)\le |H| c_H\Delta_H(u)^2
   \quad\text{for all }u.
\end{equation}
\end{enumerate}
\end{proposition}

\begin{proof}
We prove the proposition by induction on the decomposition tree.

If $H=\{e\}$ is a single edge, take $T_H=\{e\}$ and $F_H=\varnothing$. Then
\[
   Q_H(T_H;u)=-(n-1)w_e\Delta_H(u)^2,
   \quad
   Q_H(F_H;u)=w_e\Delta_H(u)^2.
\]
Thus the assertion holds with $c_H=w_e$.

Suppose next that
\[
   H=H_1\parallel\cdots\parallel H_p,
\]
and assume the result is known for all children, with constants $c_i$. Define
\[
   c_H=\frac{1}{|H|}\sum_i |H_i| c_i.
\]
For the forest certificate take
\[
   F_H=F_{H_1}\cup\cdots\cup F_{H_p}.
\]
All child terminal drops equal $\Delta_H$, so
\[
   Q_H(F_H;u)
   \le \sum_i |H_i| c_i\Delta_H^2
   =|H| c_H\Delta_H^2.
\]
For the tree certificate choose $j$ with $c_j\ge c_H$, and take
\[
   T_H=T_{H_j}\cup\bigcup_{i\ne j}F_{H_i}.
\]
Then
\begin{align*}
   Q_H(T_H;u)
   &\le -(n-|H_j|)c_j\Delta_H^2+
      \sum_{i\ne j}|H_i| c_i\Delta_H^2 \\
   &=\bigl(|H| c_H-nc_j\bigr)\Delta_H^2
   \le -(n-|H|)c_H\Delta_H^2.
\end{align*}

Finally suppose that
\[
   H=H_1\circ\cdots\circ H_s.
\]
Let $c_i$ be the child constants, and put
\[
   A_i=(n-|H_i|)c_i.
\]
Define $c_H$ by
\begin{equation}
\label{eq:c-series-combined}
   (n-|H|)c_H=\left(\sum_i\frac1{A_i}\right)^{-1}.
\end{equation}
Let $\delta_i=\Delta_{H_i}(u)$; then $\Delta_H=\sum_i\delta_i$. For the tree
certificate take $T_H=\bigcup_iT_{H_i}$. Then
\[
   Q_H(T_H;u)\le -\sum_i A_i\delta_i^2
   \le -(n-|H|)c_H\Delta_H^2
\]
by \eqref{eq:c-series-combined} and Lemma \ref{lem:harmonic-cauchy-combined}.

For a forest certificate, put
\[
   u_i=\frac1{(n-|H_i|)c_i}.
\]
Choose $j$ such that
\[
   \frac{u_j}{|H_j|}\ge\frac{\sum_i u_i}{\sum_i |H_i|}.
\]
Take
\[
   F_H=F_{H_j}\cup\bigcup_{i\ne j}T_{H_i}.
\]
With $B_j=|H_j|c_j$, the choice of $j$ is exactly the condition
\[
   \frac1{B_j}-\sum_{i\ne j}\frac1{A_i}\ge\frac1{|H| c_H}>0.
\]
Lemma \ref{lem:one-positive-combined} gives
\[
   B_j\delta_j^2-
   \sum_{i\ne j}A_i\delta_i^2
   \le |H| c_H\left(\sum_i\delta_i\right)^2,
\]
which is the desired forest inequality.
\end{proof}

We now apply the proposition at the root. Write the root parallel composition as
\[
   G=H_1\parallel\cdots\parallel H_p,
\]
and let $c_i$ be the constants supplied by Proposition
\ref{prop:tree-forest-combined} for the proper children $H_i$. Put
\[
   c_G=\frac1n\sum_i |H_i| c_i.
\]
Choose an index $j$ such that $c_j\ge c_G$, and form
\[
   \tau=T_{H_j}\cup\bigcup_{i\ne j}F_{H_i}.
\]
This is a spanning tree of $G$. The same parallel calculation as in the proof
of Proposition \ref{prop:tree-forest-combined} gives
\begin{align*}
   Q_G(\tau;u)
   &=\Energy_{E\setminus\tau}(u)-(n-1)\Energy_\tau(u)\\
   &\le -(n-|H_j|)c_j\Delta_G(u)^2+
      \sum_{i\ne j}|H_i| c_i\Delta_G(u)^2\\
   &=\bigl(nc_G-nc_j\bigr)\Delta_G(u)^2\le0.
\end{align*}
Therefore
\[
   \Energy_{E\setminus\tau}(u)\le(n-1)\Energy_\tau(u)
   \quad\text{for all }u.
\]
Together with Lemma \ref{lem:block-energy-combined}, this gives the fixed-graph
lower bound
\begin{equation}
\label{eq:fixed-graph-bound}
   \alpha_G(W)\ge\frac1n
   \quad\text{for every positive weighting }W.
\end{equation}

\subsection{Equality cases and uniqueness}

We now identify when equality can occur in \eqref{eq:fixed-graph-bound}. The
induction above naturally assigns a positive constant $c_H$ to every node of
the decomposition tree. For a leaf $H=\{e\}$, $c_H=w_e$. For a parallel node
$H=H_1\parallel\cdots\parallel H_p$,
\begin{equation}
\label{eq:c-parallel-combined}
   c_H=\frac1{|H|}\sum_i |H_i|c_{H_i}.
\end{equation}
For a series node $H=H_1\circ\cdots\circ H_s$,
\begin{equation}
\label{eq:c-series-def-combined}
   (n-|H|)c_H=
   \left(
      \sum_i\frac1{(n-|H_i|)c_{H_i}}
   \right)^{-1}.
\end{equation}

\begin{definition}
A parallel node $H=H_1\parallel\cdots\parallel H_p$ is called \emph{balanced}
if
\[
   c_{H_1}=\cdots=c_{H_p}.
\]
A series node $H=H_1\circ\cdots\circ H_s$ is called \emph{balanced} if
\[
   |H_1|(n-|H_1|)c_{H_1}
   =\cdots=
   |H_s|(n-|H_s|)c_{H_s}.
\]
The weighting is called balanced if every internal node is balanced.
\end{definition}

\begin{proposition}
\label{prop:balanced-induced-combined}
A positive weighting is balanced if and only if it is graph-induced, up to a
common scalar factor.
\end{proposition}

\begin{proof}
Let $\varphi(x)=x(n-x)$. Suppose first that the constants are balanced. Since
the root is parallel, all root children have the same constant; call it $C$.
Move down a chain
\[
   G\supset \Gamma_1\supseteq \Gamma_2\supseteq\cdots\supseteq \Gamma_{2d}=\{e\}.
\]
Across a balanced parallel node the constant does not change. Across a balanced series node $H\supset H_i$ we have
\[
   \varphi(|H_i|)c_{H_i}=\varphi(|H|)c_H,
\]
and hence
\[
   c_{H_i}=\frac{\varphi(|H|)}{\varphi(|H_i|)}c_H.
\]
Since $c_{\{e\}}=w_e$, we obtain
\[
   w_e=C
   \frac{\varphi(|\Gamma_1|)}{\varphi(|\Gamma_2|)}
   \frac{\varphi(|\Gamma_3|)}{\varphi(|\Gamma_4|)}
   \cdots
   \frac{\varphi(|\Gamma_{2d-1}|)}{\varphi(|\Gamma_{2d}|)}.
\]
This is exactly the graph-induced formula, multiplied by the common scalar
$C$.

Conversely, if the weights have this form, the same calculation read upwards
shows that constants are unchanged across parallel nodes and that
$\varphi(|H_i|)c_{H_i}$ is unchanged across the children of every series node.
Thus every internal node is balanced.
\end{proof}

\begin{lemma}\label{lem:equality-propagation-combined}
Consider any recursively chosen certificate from Proposition
\ref{prop:tree-forest-combined}, either in tree mode or in forest mode. If
equality holds in the corresponding estimate and the terminal drop of the
parent node is zero, then all child terminal drops are zero. Iterating,
equality together with $\Delta_H(u)=0$ forces all selected edge drops inside
$H$ to be zero.

If equality holds and $\Delta_H(u)\ne0$, then every child terminal drop is a
nonzero scalar multiple of $\Delta_H(u)$.
\end{lemma}

\begin{proof}
For a parallel node, all child terminal drops are equal to the parent terminal
drop, so the claim is immediate.

For a series tree certificate, equality in Cauchy's inequality occurs only when
\[
   \delta_i=\lambda/A_i
\]
for one scalar $\lambda$. Thus all child drops vanish when the parent drop
vanishes, and otherwise they are nonzero multiples of the parent drop.

For a series forest certificate, the relevant quadratic estimate is
\[
   B_j\delta_j^2-
   \sum_{i\ne j}A_i\delta_i^2
   \le
   D\left(\sum_i\delta_i\right)^2,
   \quad
   D^{-1}=\frac1{B_j}-\sum_{i\ne j}\frac1{A_i}>0.
\]
For fixed $S=\sum_i\delta_i$, the maximizer is unique. Solving the corresponding Lagrange multiplier equations gives
\[
   \delta_i=\gamma_i S
\]
with all $\gamma_i\ne0$. Hence the same conclusion holds.

The final assertion follows by induction down to the leaves. At a leaf, zero
terminal drop is exactly zero selected edge drop.
\end{proof}

\begin{lemma}\label{lem:strict-certificate-combined}
If the weighting is not balanced, then there exists a spanning tree $\tau$ and
some $\eps>0$ such that
\begin{equation}\label{eq:strict-global-combined}
   \Energy_{E\setminus\tau}(u)
   \le(n-1-\eps)\Energy_\tau(u)
   \quad\text{for every }u:V\to\R.
\end{equation}
\end{lemma}

\begin{proof}
Choose an unbalanced internal node $A$ of smallest depth. All ancestors of
$A$ are balanced.

First suppose that $A$ is parallel. We steer the recursive construction so
that $A$ occurs in tree mode: at each balanced parallel ancestor choose the
child containing $A$ as the tree child, and through a series ancestor in tree
mode all children remain in tree mode. At $A$, choose a child $A_j$ for which
$c_{A_j}=\max_i c_{A_i}$. Since $A$ is unbalanced and $c_A$ is the
$|A_i|$-weighted average of the child constants, we have $c_{A_j}>c_A$. The
parallel calculation then improves the usual tree estimate at $A$ to
\begin{equation}\label{eq:parallel-slack-combined}
   Q_A(T_A;u)
   \le -(n-|A|)c_A\Delta_A(u)^2
      -n(c_{A_j}-c_A)\Delta_A(u)^2.
\end{equation}

Now suppose that $A$ is series. Since the root is a balanced parallel node, choose at the root a
tree branch not containing $A$. The branch containing $A$ is then placed in
forest mode. At every balanced series ancestor in forest mode choose the child
containing $A$ as the forest child; at every parallel ancestor in forest mode
all children remain in forest mode. Hence $A$ occurs in forest mode. Put
\[
   u_i=\frac1{(n-|A_i|)c_{A_i}},
   \quad
   U=\sum_i u_i.
\]
Since $A$ is an unbalanced series node, the numbers $u_i/|A_i|$ are not all
equal. Choose $j$ with
\[
   \frac{u_j}{|A_j|}>\frac{U}{|A|}.
\]
In the estimate from Lemma \ref{lem:one-positive-combined} the best constant is
\[
   D_j=\left(
      \frac1{|A_j|c_{A_j}}
      -\sum_{i\ne j}\frac1{(n-|A_i|)c_{A_i}}
   \right)^{-1}.
\]
The strict inequality above is exactly $D_j<|A|c_A$. Hence the forest
estimate at $A$ improves to
\begin{equation}\label{eq:series-slack-combined}
   Q_A(F_A;u)
   \le |A|c_A\Delta_A(u)^2
      -\bigl(|A|c_A-D_j\bigr)\Delta_A(u)^2.
\end{equation}

Extend this local choice recursively to all descendants and side branches using
Proposition \ref{prop:tree-forest-combined}. At the root this produces a
spanning tree $\tau$ of $G$ and gives
\[
   Q_G(\tau;u)=\Energy_{E\setminus\tau}(u)-(n-1)\Energy_\tau(u)\le0
\]
for every $u$.

We claim that the inequality is strict for every nonconstant potential. If
not, let $u$ be nonconstant and suppose $Q_G(\tau;u)=0$. Then equality must
hold in every recursive estimate used above. In particular, equality must hold
in the strict estimate \eqref{eq:parallel-slack-combined} or
\eqref{eq:series-slack-combined}, and hence $\Delta_A(u)=0$.

By Lemma \ref{lem:equality-propagation-combined}, equality along the path from
the root to $A$ implies that every ancestor terminal drop is a nonzero multiple
of $\Delta_A(u)$. Therefore the root terminal drop is zero. Applying the same
lemma downward from the root gives zero selected edge drop on every edge of
$\tau$. Since $\tau$ is connected, $u$ is constant, a contradiction.

Thus $Q_G(\tau;u)<0$ for every nonconstant $u$. On the compact set of
potentials modulo constants normalized by $\Energy_\tau(u)=1$, the continuous
function $Q_G(\tau;u)$ attains a negative maximum; denote it by $-\eps$. This yields
\eqref{eq:strict-global-combined}.
\end{proof}

\begin{theorem}\label{thm:unique-weights-combined}
Let $G$ be a nontrivial $2$-connected series--parallel graph with $n$ edges and
$k+1$ vertices. For a positive edge weighting $W$, the following are
equivalent:
\begin{enumerate}
\item
\[
   \alpha_G(W)=
   \max_{\substack{I\subseteq E\\ |I|=k}}\lambda_{\min}P_I(W)
   =\frac1n.
\]
\item The weights are graph-induced, up to multiplication by a common positive
scalar.
\end{enumerate}
\end{theorem}

\begin{proof}
If $W$ is graph-induced, Theorem \ref{thm:main-construction} shows that
\[
   \lambda_{\min}P_\tau(W)=\frac1n
\]
for every spanning tree $\tau$. Every non-tree $k$-edge set gives a singular
block. Hence $\alpha_G(W)=1/n$.

Conversely, suppose $W$ is not graph-induced. By Proposition
\ref{prop:balanced-induced-combined}, the canonical constants are not balanced.
Lemma \ref{lem:strict-certificate-combined} gives a spanning tree $\tau$ and
$\eps>0$ such that
\[
   \Energy_{E\setminus\tau}(u)
   \le(n-1-\eps)\Energy_\tau(u)
   \quad\text{for every }u.
\]
Thus $\rho_\tau(W)\le n-1-\eps$. Lemma \ref{lem:block-energy-combined} gives
\[
   \lambda_{\min}P_\tau(W)
   =\frac1{1+\rho_\tau(W)}
   \ge\frac1{n-\eps}
   >\frac1n.
\]
Therefore $\alpha_G(W)>1/n$. This proves the contrapositive.
\end{proof}

\section{Final remarks}

The results presented here have a natural physical interpretation.

Regard $G$ as an electrical network with edge conductance $w_e$. For a voltage potential $u:V\to\R$, the Dirichlet energy
\[
   \Energy_S(u)=\sum_{xy\in S}w_{xy}\bigl(u(x)-u(y)\bigr)^2,
\]
introduced in the previous section, is the power dissipated on the set of conductors $S$.

For a spanning tree $\tau$, Lemma \ref{lem:block-energy-combined} gives
\[
   \lambda_{\min}P_\tau(W) = \inf_{u\not\equiv\mathrm{const}}
   \frac{\Energy_{\tau}(u)}{\Energy_G(u)}.
\]

Thus,
\[
   \alpha_G(W)=\max_\tau \inf_{u\not\equiv\mathrm{const}}
   \frac{\Energy_{\tau}(u)}{\Energy_G(u)}
\]
is the largest, over spanning trees, of the worst-case fraction of total energy captured by the tree.

The inequality $\alpha_G(W) \geq 1/n$ says that every positively weighted 2-sp-network contains a spanning tree capturing at least $1/n$ of the total energy for every voltage potential. The graph-induced weights are precisely the conductances for which this guarantee is sharp, while the $(G,\tau)$-induced coefficients $y_e$ are the currents through the oriented graph edges in the corresponding critical voltage mode.

Using the same terminology, the column space of transfer current matrix $Y$, also known as the star space of $G$, is the linear space of currents -- antisymmetric functions on edges -- satisfying Ohm's law and Kirchhoff's voltage law, and arising while routing arbitrary demands -- zero-sum functions on vertices.

Due to the special role of the transfer current matrix in the uniform spanning tree model (\cite{BP1993}) and more specifically, the probabilistic interpretation of the eigenvalues of its submatrices in this model (\cite{Bapat1992}), the presented results may be of interest in probability theory.

Since this study began with the discovery of the OEIS entry \cite{A115594}, we would like to mention the potential parallels between the presented results and several extremal appearances of series--parallel matroids in algebraic combinatorics, most notably Speyer's $f$-vector conjecture for tropical linear spaces \cite{SP2008}.

\section*{Acknowledgements}

The author thanks Igor Makhlin, Stanislav Budzinskiy, Tomoyuki Shirai, Aleksei Ustimenko and the authors of \cite{GTZ1997} for fruitful discussions.

Parts of Section 4 were developed with assistance from ChatGPT, model GPT-5.4 Pro.

\bibliographystyle{plain}
\bibliography{lit}

\end{document}